\documentclass{article}
\usepackage{amssymb,amsmath,amscd,amsthm}
\usepackage{mathtools}
\usepackage{ifthen}
\usepackage{pstricks}
\usepackage{pst-node}
\usepackage{graphicx}
\usepackage{psfrag}
\title{A Better Way to Deal the Cards}
\author{Mark A. Conger and Jason Howald}
\date{}
\newenvironment{piecewise}{\left\{\begin{array}{ll}}{\end{array}\right.}
\newenvironment{piecewise*}{\left\{\begin{array}{rl}}{\end{array}\right.}
\newcommand{\pwif}[2]{#1 & \mbox{if \ensuremath{#2}}}

\newcommand{\pwottt}[1]{#1 & \mbox{otherwise,}}

\newcommand{\ex}{\mathbb{E}}
\newcommand{\pr}{\mathbb{P}}

\newcommand{\abs}[1]{\ensuremath{\left| #1 \right|}}
\newcommand{\norm}[1]{\ensuremath{\left|\left| #1 \right|\right|}}
\newcommand{\orbit}[1]{\mathcal{O}(#1)}

\newcommand{\sizeof}[1]{\##1}
\newcommand{\assign}{\ensuremath{\coloneqq}}
\newcommand{\set}[1]{\left\{#1\right\}}
\newcommand{\paren}[1]{\ensuremath{\left( #1 \right)}}

\newcommand{\lit}[1]{\textbf{\texttt{#1}}}
\newcommand{\lits}[1]{\textbf{\texttt{#1}}'s}
\newcommand{\mlit}[1]{\mathbf{\mathtt{#1}}}
\newcommand{\titl}[1]{\textsl{#1}}
\newcommand{\newterm}[1]{{\bfseries #1}}

\newcommand{\euler}[2]{\genfrac{\langle}{\rangle}{0pt}{}{#1}{#2}}

\newcommand{\despoly}[2]{\mathcal{D}(#1,#2;x)}
\newcommand{\dx}{\ensuremath{\despoly{D}{D'}}}
\newcommand{\kbar}{\overline{\kappa}}
\newcommand{\des}{\mathop{\mathrm{des}}}
\newcommand{\asc}{\mathop{\mathrm{asc}}}
\newcommand{\bigo}{\mathop{\mathrm{O}}}
\newcommand{\stab}{\mathop{\mathrm{stab}}}
\newcommand{\mbar}{\overline{m}}
\newcommand{\nbar}{\overline{n}}

\newcommand{\eqnum}[1]{(\ref{eq:#1})}
\newcommand{\fignum}[1]{Figure~\ref{fig:#1}}
\newcommand{\tabnum}[1]{Table~\ref{tab:#1}}
\newcommand{\thmnum}[1]{Theorem~\ref{thm:#1}}
\newcommand{\sharpp}{\#P}

\theoremstyle{plain}
\newtheorem{theorem}{Theorem}
\theoremstyle{remark}
\newtheorem*{note}{Note}

\newcommand{\eop}{\hfill\ensuremath{\Box}} 

\newcounter{eqisnamed}
\let\realbegdispmath\[
\let\realenddispmath\]
\renewcommand{\[}[1][noname]{
	\ifthenelse{\equal{#1}{noname}}{
		\realbegdispmath
		\setcounter{eqisnamed}{0}
	}{
		\begin{equation}
		\label{eq:#1}
		\setcounter{eqisnamed}{1}
	}
}
\renewcommand{\]}{
	\ifthenelse{\equal{\value{eqisnamed}}{0}}{
		\realenddispmath
	}{
		\end{equation}
	}
}

\newcommand{\titlefont}{\rmfamily}
\newcommand{\textdivision}[4]{
	#3[#2]{\titlefont #2}
	\label{#4:#1}
}
\newcommand{\sect}[2]{\textdivision{#1}{#2}{\section}{sect}}
\newcommand{\subsect}[2]{\textdivision{#1}{#2}{\subsection}{subsect}}
\newcommand{\sectnum}[1]{Section~\ref{sect:#1}}

\newcommand{\capt}[2]{\caption[#1]{#1 #2}}

\newcommand{\psflabeltag}{label}
\newcommand{\psflabel}[3]{\psfrag{\psflabeltag#1}[#2][Bl]{#3}}
\newcommand{\ulabel}[2]{\psflabel{#1}{bc}{#2}}
\newcommand{\dlabel}[2]{\psflabel{#1}{tc}{#2}}
\newcommand{\llabel}[2]{\psflabel{#1}{cr}{#2}}
\newcommand{\rlabel}[2]{\psflabel{#1}{cl}{#2}}
\newcommand{\xlabel}[1]{\dlabel{#1}{$#1$}}
\newcommand{\ylabel}[1]{\llabel{#1}{$#1$}}

\newcommand{\tvd}{variation distance}
\newcommand{\tvdfu}{variation distance from uniform}

\newcommand{\brdeck}{\ensuremath{\mlit{N}^{13}\mlit{E}^{13}\mlit{S}^{13}\mlit{W}^{13}}}
\newcommand{\brord}{\ensuremath{D_{\mathrm{ord}}'}}
\newcommand{\brcyclic}{\ensuremath{D_{\mathrm{cyc}}'}}
\newcommand{\brrothgeb}{\ensuremath{D_{\mathrm{bf}}'}}
\newcommand{\poker}{\ensuremath{D_{\mathrm{poker}}'}}

\begin{document}
\typeout{<-- Start of my stuff -->}
\maketitle

\begin{abstract}
Most of the work on card shuffling assumes that all the cards in a deck 
are distinct, and that in a well-shuffled deck all orderings need to be 
equally likely.  We consider the case of decks with repeated cards and 
decks which are dealt into hands, as in Bridge and Poker.  We derive 
asymptotic formulas for the randomness of the resulting games.  Results 
include the influence of where a poker deck is cut, and the fact that 
switching from cyclic dealing to back-and-forth dealing will improve 
the randomness of a bridge deck by a factor of 13.
\end{abstract}

\sect{intro}{Introduction}

Card shuffling has been used as an example of a mixing problem since the 
early part of the twentieth century.  The fundamental question may be 
boiled down to: ``How fast does repeated shuffling randomize the deck?''  We 
will explore how to make mathematical sense of that question below.

At the same time, in many card games, such as bridge, euchre, and 
straight poker, the players receive ``hands'' of cards which are dealt 
from the deck after it has been shuffled.  A hand has no inherent order; 
in other words, the sequence in which the cards arrive in front of a 
player is unimportant. Thus the process of shuffling and dealing 
partitions the deck into sets of cards of predetermined sizes, and our 
goal should be to make the partition as unpredictable as possible.  This 
might not require fully randomizing the deck.

For instance, bridge is a game played with 4 players and a 52-card deck. 
Each player receives 13 cards.  The usual method is to shuffle and then 
deal cyclically: first card to the player on the dealer's left, next 
card to the player on {\it his} left, and so on, clockwise around the 
table.  Why deal this way, instead of just giving 13 cards off the top 
of the shuffled deck to the first player, then 13 to the next player, 
and so on?  If the deck were perfectly randomized by the shuffling, the 
method of dealing would not matter.

One reason not to simply ``cut the deck into hands'' is to prevent a 
dishonest dealer from stacking the deck through unfair shuffling.  But 
one might also guess, correctly, that dealing cyclically augments the 
randomness of the game when the dealer is honest but has not shuffled 
the deck thoroughly.  In this paper we address the questions: ``How much 
difference does a dealing method make?'' and ``Is cyclic dealing the 
best method?''  The goal is to show that the answers are ``Quite a bit'' 
and ``No, we can do a lot better!''  We present results for bridge and 
straight poker as examples.  We also show how to estimate the randomness 
in a deck that has distinct cards and a deck with only two types of cards.

\sect{previous}{A Brief History of Card Shuffling}

Henri Poincar\'e devoted eight sections of his 1912 book \titl{Calcul 
des Probabilit\'es} \cite[\S 225--232]{Poincare} to card shuffling.  He 
did not attempt to model any particular kind of shuffling, but showed 
that any shuffling method which meets certain mild criteria, if applied 
repeatedly to a deck, will eventually result in a well-mixed deck---that 
is, with enough shuffles the bias can be made arbitrarily small.  About 
the same time Markov \cite{Markov} was creating the more general theory 
of Markov chains, and he often used card shuffling as an example.  The 
verdict of history seems to be that Markov justly deserves credit for 
the theory named after him, but that Poincar\'e anticipated some of 
Markov's ideas in his work on card shuffling.  Most subsequent work on 
shuffling has approached the problem as a Markov chain.

In the 1950s Gilbert and Shannon \cite{Gilbert} considered the problem 
of \newterm{riffle shuffling}.  Riffle shuffling is the most common 
method used by card players to randomize a deck: the shuffler cuts the 
deck into two packets, then interleaves (riffles) them together in some 
fashion.  Using the new science of information theory, Gilbert and 
Shannon began the inquiry into how fast riffle shuffling mixes a deck. 
In the 1980s Reeds \cite{Reeds} and Aldous \cite{AldousWalks} added the 
assumption that all cut/riffle combinations are equally likely, and that 
has become known as the Gilbert-Shannon-Reeds or GSR model of card 
shuffling.

In 1992, in the most celebrated paper on card shuffling to date 
\cite{BayerDiaconis}, Bayer and Diaconis generalized the GSR shuffle to 
the \newterm{$a$-shuffle}. Let $a$ be a positive integer, and cut a deck 
into $a$ packets (one imagines an $a$-handed dealer in a futuristic 
casino), then riffle them together in some fashion.  Assuming as before 
that all cut/riffle combinations are equally likely, performing a 
randomly selected $a$-shuffle followed by a randomly selected 
$b$-shuffle turns out to be equivalent to performing a randomly selected 
$ab$-shuffle.  In particular that means that a sequence of $k$ GSR 
shuffles is equivalent to a single $2^k$-shuffle. Thus if we understand 
$a$-shuffles we implicitly understand repeated shuffles.

Bayer and Diaconis found an explicit formula for the probability of a 
particular permutation $\pi$ after an $a$-shuffle, namely
\[[probofpi] \pr_a(\pi) = \frac1{a^n} \binom{a+n-\des(\pi)-1}{n}, \]
where $n$ is the size of the deck and
\[ \des(\pi) \assign \sizeof{\set{i:\pi(i) > \pi(i+1)}} \]
is the number of descents in the permutation $\pi$.

\begin{note}
There are two ways to view a permutation:
\begin{enumerate}
\item as a bijection $\pi$ from $\set{1,2,\ldots,n}$ to iteslf, so that 
$\pi$, when applied to a sequence of objects, moves the object in 
position $i$ to position $\pi(i)$, and
\item as an ordering $\sigma = \sigma_1,\sigma_2,\ldots,\sigma_n$ of 
$1,2,\ldots,n$.
\end{enumerate}
In this paper our decks will contain repeated cards, and our 
permutations will act to rearrange them.  So we will consistently 
interpret permutations as maps.  This may disorient readers who are used 
to the other viewpoint.
\end{note}

In order to analyze the progress of a shuffler toward a well-mixed deck, 
we need a measure of how close the distribution after an $a$-shuffle is 
to the uniform distribution (all orderings equally likely).  For this 
Bayer and Diaconis use {\bf\tvdfu}, which may be defined as
\[[tvddef]
	\norm{\pr_a - U} \assign 
	\frac12 \sum_{\pi \in S_n} \abs{\pr_a(\pi) - U(\pi)},
\]
where $U$ represents the uniform distribution on permutations, i.e., 
$U(\pi) = 1/n!$ for all $\pi \in S_n$.  In terms of cards, 
the most biased game one could play with the shuffled deck is the 
two-player game in which player 1 wins whenever the ordering of the deck 
has probability higher than it should be under the uniform distribution. 
So $\norm{\pr_a - U}$ is the maximum bias toward any player in any game 
one might care to play with the shuffled deck.  That is, the probability 
of any set of permutations is within $\norm{\pr_a - U}$ of what it would 
be under the uniform distribution.

\begin{figure}
	\centering
	\rlabel{a}{$a$}
	\ulabel{y}{$y$}
	\llabel{onehalf}{$\frac12$}
	\includegraphics{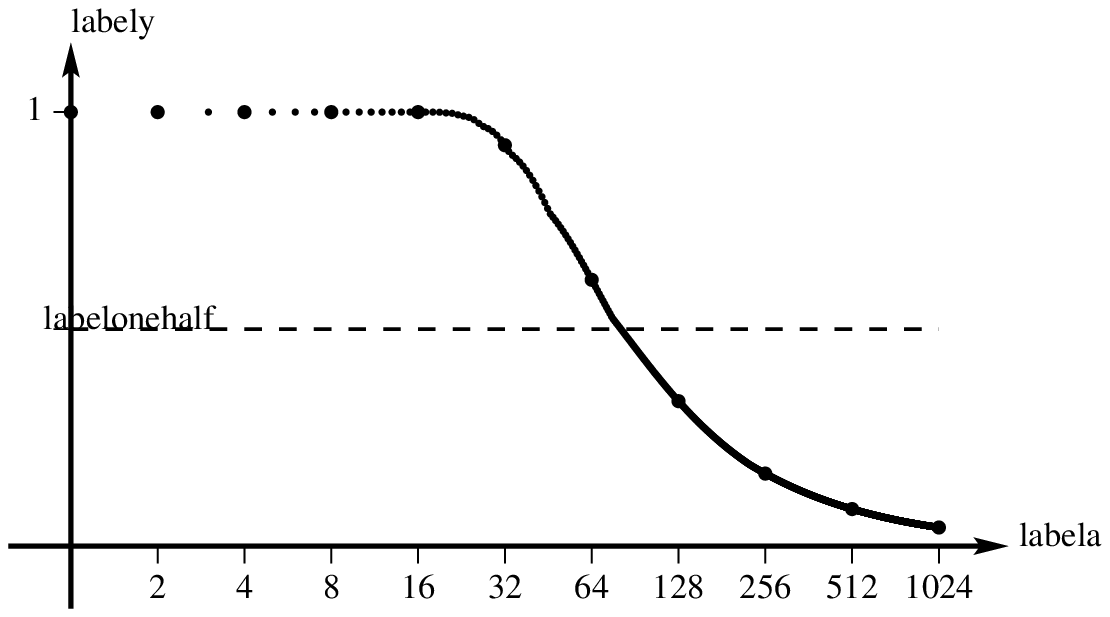}
	\caption{The \tvdfu\ of a distinct 52-card deck after an $a$-shuffle.}
	\label{fig:alldistinct}
\end{figure}

Using the probability formula in \eqnum{probofpi} and the knowledge that 
the number of permutations in $S_n$ with $d$ descents is the well-studied 
Eulerian number $\euler{n}{d}$ (see for example \cite{ConcreteMath}, 
\cite{Carlitz}, \cite{Tanny}), Bayer and Diaconis are able to compute
\[
	\norm{\pr_a - U} = \frac12
	\sum_{d=0}^{n-1}
	\euler{n}{d} \abs{
		\frac1{a^n}\binom{a+n-d-1}{n} - \frac1{n!}
	}
\]
very quickly.  The result, for $n=52$ cards and $a$ between 1 and 1024, is 
graphed in \fignum{alldistinct}.  The horizontal scale is logarithmic 
to represent the fact that a $2^k$-shuffle is the same as $k$ GSR shuffles.

The ``waterfall'' shape of the graph is typical of what Aldous 
\cite{AldousWalks} calls ``rapidly mixing Markov chains'': negligible 
change after the first few GSR shuffles, followed by a fast approach to 
uniform (``the cutoff''), and eventually halving with each extra 
shuffle.  Note that if we cast the problem in terms of $a$-shuffles, the 
ultimate exponential decay means the \tvd\ approaches $\kappa_1/a$ for 
some constant $\kappa_1$ as $a$ gets large.

Many more papers on card shuffling and its applications have been 
published.  See \cite{Diaconis01} for a survey and \cite{Mann1} for an 
excellent exposition of the Bayer and Diaconis results.  There are other 
choices besides variation distance for measuring randomness---see 
\cite[pp. 22--23]{Diaconis88}, \cite[\S5.1]{BayerDiaconis}, 
\cite[\S5]{Ciucu}, \cite[\S3]{Trefethen}, and \cite[\S1]{ADSLong} for 
alternatives.

\sect{repeatedcards}{Repeated Cards and Dealing Methods}

The implicit assumption in Bayer and Diaconis and most other work on 
card shuffling is that all cards in the deck are distinct.  (A recent 
exception is \cite{ADSLong}.) Here we consider the alternative: suppose 
a deck $D$ is a sequence of cards, each of which has a value taken from 
some fixed set of values, and we allow two cards to have the same value.  

This complicates the problem in the following ways:
\begin{itemize}
	\item Decks (ordered sequences of cards) and transformations between 
	decks can no longer be identified with permutations.  Instead for each 
	pair of decks there is a set of permutations which transform the first 
	into the second (and a different set which goes the other way). The 
	transformation sets are easy to describe, but it is difficult to find 
	their probability after shuffling.
	\item The initial order of a deck, and not just its composition, 
	affects how fast the distribution approaches uniform.  
\end{itemize}

If $D'$ is some rearrangement of $D$, let $T(D,D')$ be the set of 
permutations which, when applied to $D$, result in $D'$. ($T(D,D')$ is a 
left coset of the stabilizer $\stab(D) = T(D,D)$ and a right coset of 
$\stab(D') = T(D',D')$.) Thus the probability of obtaining $D'$ as a 
result of $a$-shuffling $D$ is
\[
	\pr_a(D \to D') \assign
	\sum_{\pi \in T(D,D')} \pr_a(\pi) =
	\frac1{a^n} \sum_d b_d \binom{a+n-d-1}{n},
\]
where $b_d$ is the number of permutations in $T(D,D')$ with $d$ descents.
In this context we call $D$ the \newterm{source deck}, $D'$ the \newterm{target 
deck}, and
\[[despolydef]
	\dx \assign
	\sum_{\pi \in T(D,D')} x^{\des(\pi)} =
	\sum_d b_d x^d
\]
the \newterm{descent polynomial} of $D$ and $D'$.  (The reader will 
kindly forgive the many different uses of the letter ``D'' in this paper. 
$D$ will always be a source deck, $D'$ a target deck, $d$ an integer 
which represents a number of descents, and $\mathcal{D}$ the descent 
polynomial.)

If we give $D$ an $a$-shuffle, then the distance between the resulting 
distribution on decks and the uniform distribution is
\[[fixedsourcetvd]
	\norm{\pr_a-U} =
	\frac12 \sum_{D' \in \orbit{D}} \abs{\pr_a(D \to D') - \frac1N},
\]
where $\orbit{D}$ is the set of reorderings of $D$ (i.e., $D$'s orbit when 
acted on by $S_n$) and $N$ is the size of $\orbit{D}$.  We will refer to 
this as the \newterm{fixed source} case.

On the other hand, suppose we are playing bridge.  Here all the cards have 
distinct values; let $e(1), e(2), \ldots, e(52)$ be the initial order of 
the cards.  The dealer $a$-shuffles the deck and then deals it out to the 
four players, who are referred to as North, East, South, and West.  Since 
each player receives 13 cards, we can describe a method of dealing cards as 
a sequence of 13 \lits{N}, 13 \lits{E}, 13 \lits{S}, and 13 \lits{W}.  For 
instance,
\[
	D' \assign
	\mlit{NESWNESW \cdots NESW} = (\lit{NESW})^{13}
\]
represents cyclic dealing, where the top card goes to North, the second 
to East, etc.  We can describe a particular partition of the deck into 
hands by a string of the same type.  For example,
\[[bridgeex]
	D \assign 
	\mlit{NSEENNWEWSSWESWNNNEESSSSSESWWNNSENWSEWSWWWEENEWNNNWE} 
\]
refers to the partition in which the North player gets cards 
$e(1),e(5),e(6),\ldots,e(50)$, the East player gets cards 
$e(3),e(4),e(8),\ldots,e(52)$, and so on.  So $D$, like $D'$, is a function 
from $\set{1,2,\ldots,52}$ to $\set{\mlit{N},\mlit{E},\mlit{S},\mlit{W}}$, 
and $D(i)$ is the player who receives card $e(i)$.

In order for the North player to receive the cards assigned by $D$, the 
shuffle must move those cards to the positions occupied by \lits{N} in 
$D'$, and likewise for the other players. In other words, if we think of 
$D$ and $D'$ as \emph{decks} with cards of value \lit{N}, \lit{E}, \lit{S}, 
and \lit{W}, the shuffle will produce the desired partition if and only if 
it takes $D$ to $D'$. Thus we already have a notation for the distribution 
over partitions that an $a$-shuffle followed by the dealing method $D'$ 
produces, and the \tvd\ between that distribution and uniform is
\[[fixedtargettvd]
	\norm{\pr_a-U} =
	\frac12 \sum_{D \in \orbit{D'}} \abs{\pr_a(D \to D') - \frac1N},
\]
where $N$ is now the number of possible partitions.  This is the 
\newterm{fixed target} case.  Note that despite the strong similarity 
between \eqnum{fixedsourcetvd} and \eqnum{fixedtargettvd}, fixed source 
and fixed target are dual problems, not identical, because the 
transition probability $\pr_a(D \to D')$ is not symmetric.

\sect{transprob}{Calculating Transition Probabilities}

Unfortunately, computing $\norm{\pr_a-U}$ precisely for 
realistically sized decks is prohibitively complicated in both the fixed 
source and fixed target cases. Both cases require knowledge of the 
transition probabilities $\pr_a(D \to D')$ in order to calculate \tvdfu. 
Conger and Viswanath \cite{Shuffle1} showed that calculation of 
transition probabilities is computationally equivalent to calculation of 
the coefficients of the descent polynomial $\dx$.  Unfortunately, the 
same authors subsequently found \cite{Shuffle3} that for certain decks 
the calculation belongs to a class of counting problems called \sharpp, 
and is in fact \sharpp-complete.  As with  NP-complete problems, it is 
generally believed that \sharpp-complete problems do not admit efficient 
solutions.

Barring a method for calculating variation distance without first 
computing transition probabilities, the question shifts to 
approximation.  \thmnum{WZ} (below) will allow us to approximate 
transition probabilities when $a$ is large, given some simple 
information about the decks.  To that end, here are some definitions:

If $u$ and $v$ are card values, we say that $D$ has a \newterm{$u$-$v$ 
digraph at $i$} if $D(i)=u$ and $D(i+1)=v$.  We say that $D$ has a 
\newterm{$u$-$v$ pair at $(i,j)$} if $i<j$, $D(i)=u$, and $D(j)=v$.  
The distinction between digraphs and pairs is akin to that between descents 
and inversions in a permutation.  Let
\begin{align}
	W(D,u,v) &\assign \label{eq:wdef}
	\sizeof{\set{\mbox{$u$-$v$ digraphs in $D$}}} -
	\sizeof{\set{\mbox{$v$-$u$ digraphs in $D$}}}, \\
	Z(D,u,v) &\assign \label{eq:zdef}
	\sizeof{\set{\mbox{$u$-$v$ pairs in $D$}}} -
	\sizeof{\set{\mbox{$v$-$u$ pairs in $D$}}}.
\end{align}
For example, the deck $D = \mlit{ABAAABABB}$ has 3 \lit{A-B} digraphs and 2
\lit{B-A} digraphs, so $W(D,\mlit{A},\mlit{B}) = 3-2 = 1$.  $D$ has 15
\lit{A-B} pairs and 5 \lit{B-A} pairs, so $Z(D,\mlit{A},\mlit{B}) = 15-5 =
10$.  Note that both $W$ and $Z$ are antisymmetric in $u$ and $v$:
$W(D,u,v) = -W(D,v,u)$ and $Z(D,u,v) = -Z(D,v,u)$.

In the theorem below we assume that $D$ is a deck of $n$ cards.  For each 
card value $v$, $n_v$ is the number of cards in $D$ with that value.  For 
convenience we assume an implicit order on the values; the particular order 
chosen is arbitrary and does not affect the result.

\begin{theorem} \label{thm:WZ}
	Suppose $D$ is as above, and $D'$ is a reordering of $D$.  Then
	\[
		\pr_a(D \to D') =
		\frac1N + c_1(D,D')a^{-1} + \bigo(a^{-2}),
	\]
	where $N$ is the number of reorderings of $D$ and
	\[[c1def]
		c_1(D,D') =
		\frac{n}{2N}
		\sum_{u<v}
		\frac{W(D,u,v)Z(D',u,v)}{n_un_v}.
	\]
\end{theorem}

We begin with a plausibility argument in favor of formula \eqnum{c1def}. 
An $a$-shuffle begins with an $a$-cut, which arranges the cards into $a$ 
piles.  If $a$ is very large, most piles will have size $0$ or $1$.  If 
no pile has two or more cards, then no ordering information will survive 
the shuffle.  The next most likely case is that the $a$-cut produces 
just one pile with two cards $u$ and $v$.  Should this happen, these cards 
must come from a $u$-$v$ digraph in $D$. When the cards are reassembled, 
these two must remain in order, becoming a $u$-$v$ pair in $D'$. Thus 
the main source of bias in $\pr_a(D \to D')$ is the relationship between 
digraphs in $D$ and pairs in $D'$. This suggests the formula 
$W(D,u,v)Z(D',u,v)$.

\emph{Proof of \thmnum{WZ}.} Let $S$ be the size of the stabilizer of 
$D$ in $S_n$.  Since $T(D,D')$ is a coset of the stabilizer, its size is 
$S$ also, and since there are $N$ such cosets, $NS = n!$.  Let $b_d$ be 
the number of permutations in $T(D,D')$ with $d$ descents; so $\sum_d 
b_d = S$.  Then
\begin{align*}
	\pr_a(D \to D') &=
	\frac1{a^n} \sum_d b_d \binom{n+a-d-1}{n} \\&=
	\frac1{a^n} \sum_d b_d \frac{(a-d)(a-d+1)\cdots(a-d+n-1)}{n!} \\&=
	\frac1{n!} \sum_d b_d
		\paren{1 + \frac{-d}{a}}
		\paren{1 + \frac{1-d}{a}} \cdots
		\paren{1 + \frac{n-1-d}{a}} \\&=
	\frac1{NS} \paren{
		\sum_d b_d + a^{-1} \sum_d b_d \paren{
			\frac{n(n-1)}2 - nd
		} + \bigo(a^{-2})
	}.
\end{align*}
So the constant term is $1/N$, which is what we expect: it means that 
if we shuffle $D$ for long enough, the probability of obtaining any 
particular deck approaches $1/N$, i.e., the distribution on decks 
approaches uniform.  The coefficient of $a^{-1}$ is
\[
	\frac{n}{2NS} \sum_d b_d (n-1-2d) =
	\frac{n}{2N} \ex(n-1-2\des(\pi)),
\]
where $\pi$ is a permutation chosen uniformly from $T(D,D')$ and $\ex$ 
represents expectation.  Recall that descents are positions $i$ with $1 \le 
i \le n-1$ such that $\pi(i) > \pi(i+1)$.  The other positions are 
\newterm{ascents}, and if we denote their number by $\asc(\pi)$ then we 
must have $\des(\pi) + \asc(\pi) = n-1$.  So the number we are after is
\[[omegasum]
	\frac{n}{2N} \ex(\asc(\pi)-\des(\pi)) =
	\frac{n}{2N} \sum_{i=1}^{n-1} \ex{\omega_i(\pi)},
\]
where
\[
	\omega_i(\pi) = \begin{piecewise*}
		\pwif{1}{\pi(i) < \pi(i+1),} \\
		\pwif{-1}{\pi(i) > \pi(i+1).}
	\end{piecewise*}
\]
Let the first card in $D$ have value $u$ and the second have value $v$.  
Suppose first that $u=v$.  $\pi$ must take those two cards to two positions 
in $D'$ which have value $u$, but otherwise it has no reason to prefer any 
particular destinations; thus it is equally likely that $\pi(1)<\pi(2)$ as 
that $\pi(1)>\pi(2)$.  So if $u=v$ then $\ex \omega_1(\pi) = 0$.

On the other hand, suppose $u \ne v$.  Then $\pi$ picks a destination for 
the top card uniformly from among those $j$ for which $D'(j)=u$, and 
likewise $\pi(2)$ is chosen uniformly and independently from $D'^{-1}(v)$. 
We will have $\omega_1(\pi)=1$ if the first choice is less than the second, 
that is, if $\pi$ maps $\set{1,2}$ to a $u$-$v$ pair in $D'$. 
$\omega_1(\pi)$ will be $-1$ if $\pi$ maps $\set{1,2}$ to a $v$-$u$ pair in 
$D'$.  Each pair is equally likely, so
\[
	\ex \omega_1(\pi) = \frac{
		\sizeof{\set{\mbox{$u$-$v$ pairs in $D'$}}} -
		\sizeof{\set{\mbox{$v$-$u$ pairs in $D'$}}}
	}{
		\sizeof{\set{\mbox{$u$-$v$ pairs in $D'$}}} +
		\sizeof{\set{\mbox{$v$-$u$ pairs in $D'$}}}
	} =
	\frac{Z(D',u,v)}{n_un_v}.
\]
All the other $\omega_i$ are calculated in the same way, and if we group 
them according to the values of the digraphs in $D$, we get the desired 
result for $c_1$. \eop

\sect{kappa1}{The Fixed Source Case}

Now that we can approximate the transition probability between two 
decks, we can approximate variation distances. In the case of a fixed 
source deck $D$ we have
\[
	\norm{\pr_a - U} =
	\frac12 \sum_{D' \in \orbit{D}}
	\abs{\pr_a(D \to D') - \frac1N} =
	\kappa_1 a^{-1} + \bigo(a^{-2}),
\]
where
\[[kappa1def]
	\kappa_1 =
	\kappa_1(D) \assign
	\frac12 \sum_{D' \in \orbit{D}} \abs{c_1(D,D')}.
\]

\subsect{alldistinct}{All-Distinct Decks}

\begin{figure}
	\centering
	\rlabel{a}{$a$}
	\ulabel{y}{$y$}
	\xlabel{1}
	\xlabel{2}
	\xlabel{4}
	\xlabel{8}
	\xlabel{16}
	\xlabel{32}
	\xlabel{64}
	\xlabel{128}
	\xlabel{256}
	\xlabel{512}
	\xlabel{1024}
	\ylabel{1}
	\llabel{onehalf}{$\frac12$}
	\includegraphics{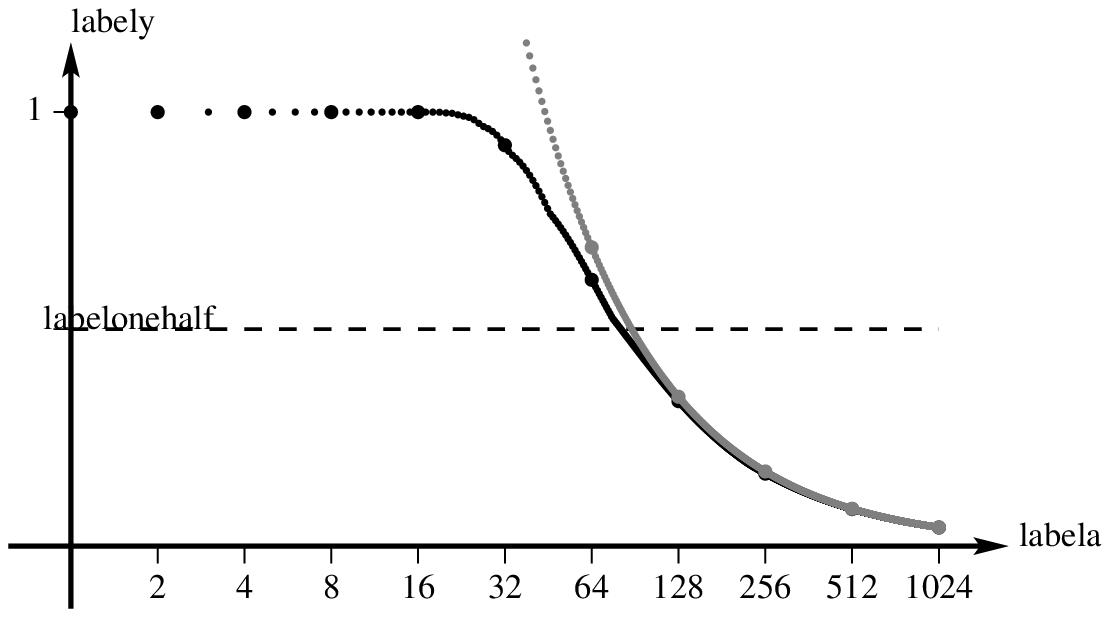}
	\capt{The \tvd\ from uniform, and first-order approximation, of a 
	distinct 52-card deck after an $a$-shuffle.} {The actual \tvd\ is 
	graphed in black, and the first order approximation $\kappa_1/a$, with 
	$\kappa_1=44.05710497$, is graphed in gray.}
	\label{fig:alldistinctwithkappa1}
\end{figure}

For example, consider shuffling a deck $D$ containing $n$ distinct cards. 
Without loss of generality we may assume that $D = 1,2,\ldots,n$, i.e., 
that the value of the card in position $i$ is $i$. Each reordering of $D$ 
is produced by a unique permutation, so
\[ \kappa_1(D) = \frac12 \sum_{\pi \in S_n} \abs{c_1(D, \pi D)}. \]
Since each card appears once and all the digraphs in $D$ are of the form 
$(i,i+1)$ we can reduce \eqnum{c1def} to
\[
	c_1(D,\pi D) = \frac{n}{2N} \sum_{i=1}^{n-1} Z(\pi D, i, i+1).
\]
$Z(\pi D,i,i+1)$ is 1 if $\pi$ has an ascent at $i$ and $-1$ if $\pi$ has a 
descent at $i$, so we have
\[[alldistinctk1]
	\kappa_1(D) =
	\frac{n}{4N} \sum_{\pi \in S_n} \abs{\asc(\pi)-\des(\pi)} =
	\frac{n}{4N} \sum_d \euler{n}{d} \abs{n-1-2d}.
\]
The Eulerian numbers can be calculated using a simple recurrence 
\cite{ConcreteMath}, so \eqnum{alldistinctk1} is all we need to find the 
long-term behavior of \tvdfu\ after shuffling a deck of $n$ distinct cards.  
For instance, when $n=52$ we find that $\kappa_1$ is
\[
\frac{146020943891326775423340146124729913263177343486982212261189487693}
{3314356310443124530393681659122442758682178888925184000000000000},
\]
which is approximately 44.06.  (In  general, $\kappa_1$ is about
$n\sqrt{(n+1)/24\pi}$ for a deck of $n$ distinct cards \cite{Thesis}.)
So after giving a deck of 52 distinct cards 
an $a$-shuffle, where $a$ is large, the \tvd\ from uniform will be 
approximately $44.06/a$.  We can compare that with the exact results Bayer 
and Diaconis obtained for the same deck, to see how big $a$ has to be to 
make the approximation a good one.  \fignum{alldistinctwithkappa1} shows 
the result.

\subsect{twovalues}{Decks with Two Card Types}

Let $G_k$ be a complete directed graph on $k$ vertices, with loops at 
each vertex.  A deck with cards of $k$ different types may be thought of 
as a walk on $G_k$, where the starting position is the top card and the 
ending position is the bottom card.  Each edge in the walk represents a 
digraph in the deck.

\begin{center}
\psset{unit=1pt}
\begin{pspicture}(0,0)(150,50)
	\rput{0}(25,25){\circlenode{one}{\Large 1}}
	\rput{0}(125,25){\circlenode{two}{\Large 2}}
	\nccircle[angleA=90,nodesepA=3]{<-}{one}{15}
	\ncarc[arcangleA=30,arcangleB=30,nodesepA=3,nodesepB=3]{->}{one}{two}
	\ncarc[arcangleA=30,arcangleB=30,nodesepA=3,nodesepB=3]{->}{two}{one}
	\nccircle[angleA=270,nodesepA=3]{<-}{two}{15}
\end{pspicture}
\end{center}

Consider a deck $D$ with only two types of cards, which we will 
unimaginatively label 1 and 2.  If the deck begins with 1 and ends with 
2, then the corresponding walk must have traversed the 1-2 edge once more 
than the 2-1 edge, and therefore $W(D,1,2)=1$.  Likewise beginning with a 2 
and ending with a 1 makes $W(D,1,2)=-1$.  But beginning and ending with the 
same type of card means that both edges were traversed the same number of 
times, so in that case $W(D,1,2)$ is 0.  With only two types of cards 
\eqnum{c1def} reduces to
\[
	c_1(D,D') =
	\frac{n}{2n_1n_2N} W(D,1,2) Z(D',1,2),
\]
which vanishes if $W(D,1,2)$ is 0, for all $D'$.  So if the top and 
bottom cards of the unshuffled deck are the same, $\kappa_1$ will be 0, 
meaning that in the long run the variation distance decreases at least as 
fast as some multiple of $a^{-2}$.  So we have the surprising result 
that a shuffler of this type of deck can greatly speed the deck's 
approach to randomness by making sure that the top and bottom cards are 
the same before he begins shuffling.

We will concentrate on the fixed target case for the rest of this paper, 
but the interested reader may consult \cite{Thesis} for more results in 
the fixed source case.

\sect{fixedtarget}{The Fixed Target Case (Dealing into Hands)}

If we fix the target deck (i.e., dealing method) $D'$, we have

\[
	\norm{\pr_a - U} =
	\frac12\sum_{D \in \orbit{D'}} \abs{\pr_a(D \to D') - \frac1N} =
	\kbar_1 a^{-1} + \bigo(a^{-2}),
\]
where
\[
	\kbar_1 =
	\kbar_1(D') :=
	\frac12\sum_{D \in \orbit{D'}} \abs{c_1(D,D')}.
\]
This sum seems intractably large, so it is useful to have an alternative algorithm 
for the calculation of $\kbar_1(D')$.  First notice that 
\begin{align*}
	\kbar_1(D') &=
	\frac12 \sum_{D \in \orbit{D'}} \abs{c_1(D,D')} \\&=
	\frac12 \sum_{D \in \orbit{D'}} \abs{
		\frac{n}{2N} \sum_{u<v} \frac{W(D,u,v)Z(D',u,v)}{n_un_v}
	} \\&=
	\frac{n}{4N} \sum_{D \in \orbit{D'}} \abs{
		\sum_i \frac{Z(D',D(i),D(i+1))}{n_{D(i)}n_{D(i+1)}}
	} \\&=
	\frac{n}{4N} \sum_{D \in \orbit{D'}} \abs{\theta(D)},
\end{align*}
where
\[
	\theta(D) \assign \sum_{i} \frac{Z(D',D(i),D(i+1))}{n_{D(i)}n_{D(i+1)}}.
\]
For our cases, there are far fewer possible values for $\theta(D)$ than 
there are decks $D \in \orbit{D'}$, so we want to reason about the 
distribution of values of $\theta(D)$ as $D$ ranges over $\orbit{D'}$. 
To prepare a recursion, we will need this distribution to depend also on 
the last card of $D$, and we will need to consider decks with fewer 
cards than $D$.  Since $\orbit{D'}$ depends on the number of cards of 
each type in $D'$, but not on their order, let $\mbar = (m_1,m_2, 
\ldots, m_k)$ be an integer vector, representing a collection of $m_1 \le n_1$ 
cards labeled $1$, $m_2 \le n_2$ cards labeled $2$, etc.  If $v \in \{1, 
\ldots, k\}$ is a card value, we write $D \dashv v$ to mean that the 
last card of $D$ is $v$, and $D \dashv uv$ to mean that $D$ ends with 
the digraph $uv$.  We also write:
\begin{align*}
	\orbit{\mbar} &\assign \begin{piecewise}
		\pwif{\orbit{1^{m_1}2^{m_2}\cdots k^{m_k}}}{\mbox{all $m_i \ge 0$,}} \\
		\pwottt{\emptyset}
	\end{piecewise} \\
	g_{\mbar,v}(t) &\assign \sum_{
		\substack{
			D \in \orbit{\mbar} \\ D \dashv v
		}
	} t^{\theta(D)}.
\end{align*}
Then $\sum_v g_{\nbar,v}(t)$ will record the distribution 
of interest.

Let $e_v$ be the standard basis vector with a 1 in coordinate $v$ and 0s 
in all other coordinates. We derive the recurrence by considering the 
second-to-last card:
\begin{align} \label{eq:kbthetarec}
	g_{\mbar,v}(t) &=
	\sum_{
		\substack{
			D \in \orbit{\mbar} \\
			D \dashv v
		}
	} t^{\theta(D)} =
	\sum_u \sum_{
		\substack{
			D \in \orbit{\mbar} \\
			D \dashv uv
		}
	} t^{\theta(D)} \\&= \nonumber
	\sum_u \sum_{
		\substack{
			\tilde{D} \in \orbit{\mbar - e_v} \\
			\tilde{D} \dashv u
		}
	} t^{\theta(\tilde{D})+\frac{Z(D',u,v)}{n_un_v}} =
	\sum_u t^{\frac{Z(D',u,v)}{n_un_v}} g_{\mbar - e_v,u}.
\end{align}

This enables us to find $\kbar_1(D')$ by recursively 
computing $g_{\mbar,v}(t)$ for each card value $v$ and for each 
integer vector $\mbar$ with $(0,\ldots,0) \leq \mbar \leq 
(n_1, \ldots, n_k)$.  There are $k\prod (n_i+1)$ generating functions to 
compute, which is feasible for the cases we discuss.

\subsect{poker}{Straight Poker}

Straight poker is a game in which players receive 5 cards each from a deck 
of 52 distinct cards.  The remaining cards are unused.  Dealing is 
traditionally cyclic, so with 4 players the normal deal sequence is
\[ \poker = (1234)^55^{32}, \]
where 1, 2, 3, and 4 represent the players and 5 is the ``hand'' of 
unused cards.  The reader may check that $Z(\poker,u,v)$ is the row $u$, 
column $v$ entry of
\[
	Z(\poker) = \paren{
		\begin{array}{rrrrr}
			0 & 5 & 5 & 5 & \phantom{-}160 \\
			-5 & 0 & 5 & 5 & 160 \\
			-5 & -5 & 0 & 5 & 160 \\
			-5 & -5 & -5 & 0 & 160 \\
			-160 & -160 & -160 & -160 & 0
		\end{array}
	},
\]
which allows us to use \eqnum{kbthetarec} to calculate
\[
	\kbar_1(\poker) =
	\frac{1041539930128654272599}{123600572196960202344} \approx
	8.427.
\]
The usual procedure in poker, however, is for the dealer to shuffle the 
cards and then allow the player to his left to cut them---that is, move 
some number of cards from the top to the bottom---before dealing.  
Fulman \cite{FulmanCuts} showed that with a deck of distinct cards, a 
shuffle followed by a random cut was no more effective a randomizer than 
the shuffle alone.

\begin{figure}
	\centering
	\rlabel{m}{$m$}
	\ulabel{k}{$\kbar_1$}
	\includegraphics{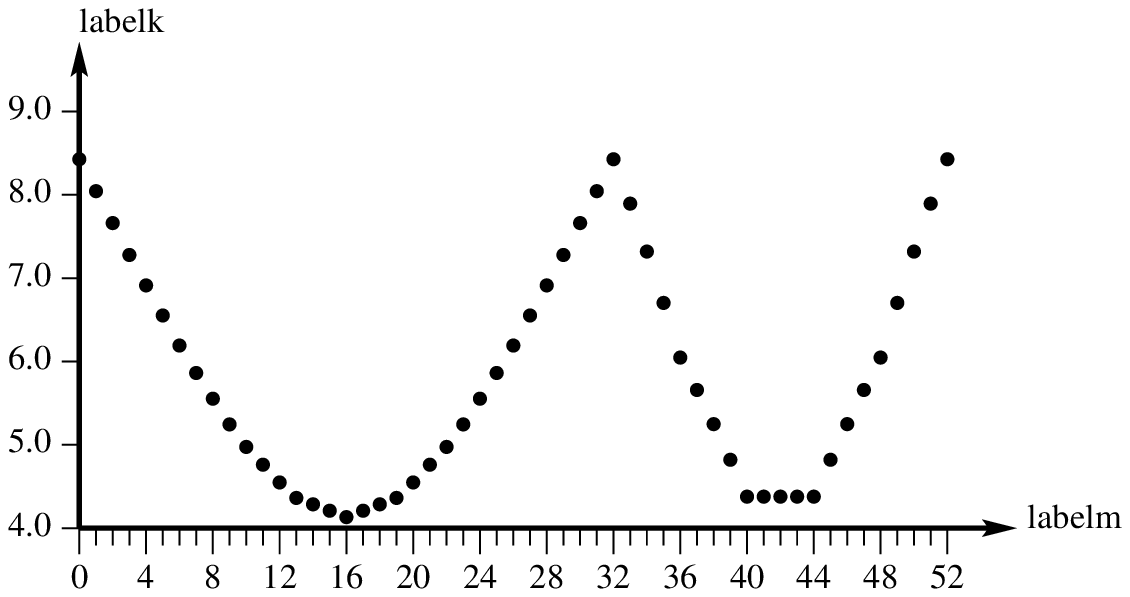}
	\capt{The effect on $\kbar_1$ of cutting a poker deck at position 
	$m$.}{}
	\label{fig:kappa1pokercuts}
\end{figure}

When the deck will be dealt into hands, however, choosing a particular 
cut \emph{can} enhance the randomness. Moving $k$ cards has the same 
effect as making the deal sequence $\sigma^k \poker$, where $\sigma$ is 
the cycle $(1,2,\ldots,52)$.  The problem is small enough that we can 
simply try all possible cuts and report $\kbar_1$ for each; the result 
is in \fignum{kappa1pokercuts}. The best place to cut the deck is after 
the 16th card, making $Z(u,5)=0$ for each player $u$.  We find \[ 
\kbar_1(\sigma^{16} D_{\mathrm{poker}}') = 
\frac{523485619699747366033}{126685078454994859800} \approx 4.132. \] 
Thus a good cut can effectively halve the value of $\kbar_1$, meaning it 
is worth one extra GSR shuffle.  The method of the next section can be 
used to improve the situation further.

\subsect{bridge}{Bridge}

\newcommand{\brmatrix}{
	\paren{
		\begin{array}{rrr@{\hspace{19pt}}r}
			0 & 1 & 1 & 1 \\
			-1 & 0 & 1 & 1 \\
			-1 & -1 & 0 & 1 \\
			-1 & -1 & -1 & 0
		\end{array}
	}
}
As described in \sectnum{repeatedcards}, a method of dealing bridge can be 
identified with a target deck $D' \in \orbit{\brdeck}$.  We have $n=52$ and 
$n_v=13$ for each card value $v$, so
\[[bridgekappa]
	\kbar_1(D') =
	\frac1{13N} \sum_{D \in \orbit{D'}}
	\abs{\sum_i Z(D',D(i),D(i+1))}.
\]

Suppose we deal a game of brige by ``cutting the deck into hands.''  That 
is, the top 13 cards go to North, the next 13 to East, etc.  Call this 
``ordered dealing.''  Symbolically,
\[ \brord = \brdeck. \]
There are 169 \lit{N}-\lit{E} pairs in $\brord$ and no \lit{E}-\lit{N} 
pairs, so $Z(\brord,\mlit{N},\mlit{E}) = 169$.  Likewise for the other card 
values, so
\[ Z(\brord) = 169 \brmatrix, \]
where we give the values the implicit ordering $\mlit{N} < \mlit{E} < 
\mlit{S} < \mlit{W}$ and interpret the entries in the matrix accordingly.  
Using \eqnum{kbthetarec} we can compute
\[
	\kbar_1(\brord) =
	\frac{93574839271687495932003418573}{3352796110343049552452340000} \approx
	27.91.
\]
Thus in the long run only slightly less shuffling is required for a bridge 
deck that will be cut into hands than for a deck in which all orderings are 
distinct.

Of course the way that most bridge players deal is cyclically:
\[ D' = \brcyclic = (\mlit{NESW})^{13}. \]
In that case the reader can check that
\[
	Z(\brcyclic) = 13 \brmatrix,
\]
which is to say, $Z(\brcyclic,u,v) = \frac1{13} Z(\brord,u,v)$ for all
pairs of card types.  It follows then from \eqnum{bridgekappa} that
\[
	\kbar_1(\brcyclic) =
	\frac1{13} \kbar_1(\brord) =
	\frac{7198064559360576610154109121}{3352796110343049552452340000} \approx
	2.147.
\]
So dealing cyclically works 13 times as well as simply cutting the deck 
into hands.  That is to say, in the long run, switching from ordered to 
cyclic dealing is worth an extra $\log_2(13) \approx 3.7$ 2-shuffles.

The reason cyclic dealing is so much better than cutting into hands is 
that it makes $Z(D,u,v)$ small by better balancing the number of $u$-$v$ 
pairs with the number of $v$-$u$ pairs, for all $u$ and $v$.  If $a$ is 
considerably larger than the deck size $n$, it is likely that when the 
deck is partitioned into $a$ packets, all of the packets will either be 
empty or contain exactly one card.  If such is the case, they will be 
riffled together in an arbitrary order, and the deck will be perfectly 
randomized.  In fact, that suggests an alternate way to estimate how 
many shuffles are needed to adequately randomize a deck, an idea due to 
Reeds \cite{Reeds} and reported in Diaconis \cite{Diaconis88}. Simply 
calculate the likelihood that each card is in a different packet after 
the cut (this is the celebrated ``birthday problem'' of combinatorics), 
and pick $a$ large enough that the probability is high.

The new idea here is that a dealing method can help ameliorate the bias in 
the case where $a$ is still reasonably large, but some packet contains two 
cards.  Imagine that the cards are initially arranged from ``best'' to 
``worst'' before shuffling.  Then a packet with two cards contains a 
``good'' card atop a ``worse'' one, and after riffling the two cards will 
remain in the same order, though other cards may come between them.  This 
is the source of the bias which remains even when $a$ is large.  If the two 
cards are dealt to players $u$ and $v$, then we would like it to be 
approximately equally likely that $u$ gets the good card and $v$ the bad as 
the other way around.  Thus we would like there to be about as many $u$-$v$ 
pairs in the dealing method as there are $v$-$u$ pairs.

Consider just the North and East players in bridge.  We can describe a 
dealing method $D'$ as it applies to those players by drawing a north-east 
lattice path starting from the lower-left corner of a $13 \times 13$ grid. 
That is, traverse $D'$ and draw a north segment whenever an \lit{N} is 
encountered and an east segment whenever an \lit{E} is encountered.  (The 
traditional names of players are very fortuitous for this exercise.) Every 
square to the southeast of the path has a northward segment to its left and 
an eastward segment above it, so it corresponds to a \lit{N}-\lit{E} pair 
in $D'$.  Likewise the squares in the Young shape to the northwest of the 
path represent \lit{E}-\lit{N} pairs.

\begin{figure}
	\centering
	\includegraphics{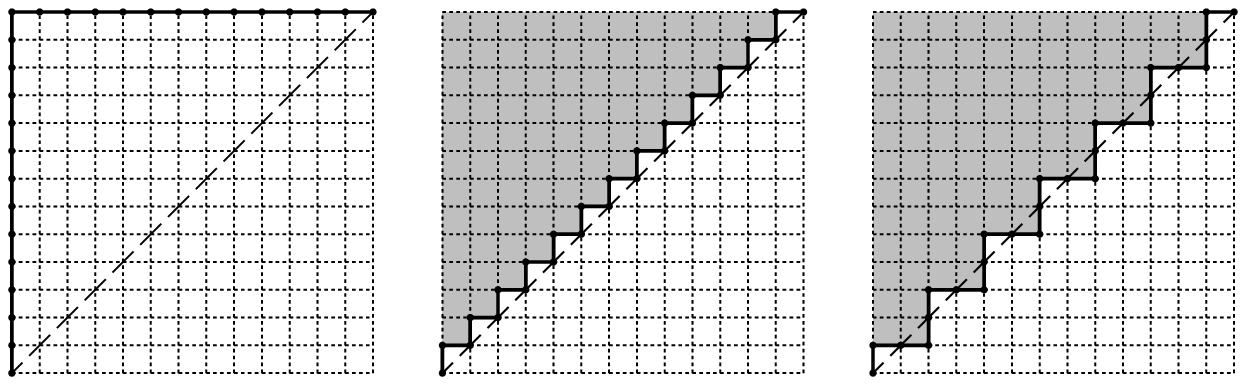}
	\capt{Three styles of bridge dealing, represented by lattice 
	paths.}{The grid on the left represents ordered dealing, the center is 
	cyclic dealing, and the one on the right is back-and-forth dealing. 
	Each grid shows the sequence of \lit{N} and \lit{E} cards as north and 
	east line segments respectively.  The size of the Young shape to the 
	northwest of the path is the number of \lit{E}-\lit{N} pairs in the 
	target deck, and the size of the complementary shape is the number of 
	\lit{N}-\lit{E} pairs.}
	\label{fig:bridgeyoung}
\end{figure}

\fignum{bridgeyoung} shows the paths and shapes for ordered and cyclic 
dealing on the left and in the middle.  Cyclic dealing is much better than 
ordered dealing because the path stays near the diagonal of the grid, so
about half the squares are on either side.

However, it always stays \emph{to one side} of the diagonal, and thus it 
is easy to see that we can do better!  The path on the right side of 
\fignum{bridgeyoung} corresponds to
\[ \brrothgeb = (\mlit{NESWWSEN})^6 \mlit{NESW}, \]
and by crossing the diagonal it balances the two sides as well as can be 
done, making $Z(\brrothgeb,\mlit{N},\mlit{E}) = 1$.  Likewise for the other 
pairs of players, so we have
\[ Z(\brrothgeb) = \brmatrix \]
and therefore
\[
	\kbar_1(\brrothgeb) =
	\frac1{13} \kbar_1(\brcyclic) =
	\frac{7198064559360576610154109121}{43586349434459644181880420000} \approx
	0.165.
\]
The dealing method described by $\brrothgeb$ may be called 
``back-and-forth'' dealing, since the dealer hands out cards once around 
the table clockwise, then once counterclockwise, then clockwise, 
counterclockwise, etc.  We have shown that in the long run 
(i.e., for large $a$), back-and-forth dealing is 13 times as effective as 
cyclic dealing.  Or, switching from cyclic to back-and-forth dealing is 
worth $\log_2(13) \approx 3.7$ extra GSR shuffles.

We could apply the same strategy to the poker deck that was cut after the 
16th card; the combination of that cut and back-and-forth dealing produces 
a $\kbar_1$ which is $1/5$ what it was with cyclic dealing.

\sect{questions}{Two Big Questions}

\begin{enumerate} \item {\bf How big does $a$ have to be for the 
first-order approximation to be a good one?}

By manipulating absolute value signs, one can show \cite{Thesis} that the 
error between \tvdfu\ and the first order estimate $\kappa_1 a^{-1}$ or 
$\kbar_1 a^{-1}$ is bounded above by
\[
	\frac12 \sum_d \euler{n}{d} \abs{
		\frac1{a^n} \binom{a+n-d-1}{n} -
		\frac1{n!} -
		\frac{1}{a(n-1)!}\paren{\frac{n-1}2 - d}
	}.
\]
This bound, which depends only on the size of the deck ($n$) and the size 
of the shuffle ($a$), could undoubtedly be improved.  
\fignum{alldistinctwithkappa1} shows that in the case of a deck of 52 
distinct cards, the first-order estimate becomes quite good at about the 
point of the cutoff.

So, how many times does a bridge player need to shuffle in order to see 
the promised benefits of switching dealing methods?  Monte Carlo 
estimates (see \cite{Thesis} for a full explanation of methods and 
confidences) give strong evidence that back-and-forth dealing beats 
cyclic dealing after 5 or more shuffles.  \tabnum{ibridge} shows the 
results.

\begin{table} \label{tab:ibridge}
	\centering
	\footnotesize
	\addtolength{\tabcolsep}{-2pt}
	\begin{tabular}{|l|l|c|c|c|c|c|c|c|} \hline
		Deck & Method &
		$a=16$ & 32 & 64 & 128 & 256 & 512 & 1024 \\ \hline
		52 Distinct & Exact & 
		1.0000 & 0.9237 & 0.6135 & 0.3341 & 0.1672 & 0.0854 & 0.0429 \\
		$123\cdots(52)$ & $44.0571a^{-1}$ &
		2.7536 & 1.3768 & 0.6884 & 0.3442 & 0.1721 & 0.0860 & 0.0430 \\ \hline
		Ordered Bridge & Monte Carlo &
		0.9902 & 0.7477 & 0.4230 & 0.2183 & 0.1104 & 0.0550 & 0.0274 \\
		$\mlit{N}^{13}\mlit{E}^{13}\mlit{S}^{13}\mlit{W}^{13}$ & $27.9095a^{-1}$ &
		1.7443 & 0.8722 & 0.4361 & 0.2180 & 0.1090 & 0.0545 & 0.0273 \\ \hline
		Cyclic Bridge & Monte Carlo & 
		0.2349 & 0.0735 & 0.0346 & 0.0169 & 0.0084 & 0.0042 & 0.0021 \\
		$(\mlit{NESW})^{13}$ & $2.1469a^{-1}$ &
		0.1342 & 0.0671 & 0.0335 & 0.0168 & 0.0084 & 0.0042 & 0.0021 \\ \hline
		Back-Forth Bridge & Monte Carlo &
		0.3118 & 0.0260 & 0.0073 & 0.0022 & 0.0008 & 0.0003 & 0.0002 \\
		$(\mlit{NESWWSEN})^6(\mlit{NESW})$ & $0.1651a^{-1}$ &
		0.0103 & 0.0052 & 0.0026 & 0.0013 & 0.0006 & 0.0003 & 0.0002 \\ \hline
	\end{tabular}
	\addtolength{\tabcolsep}{1pt}
	\caption{Variation distances from uniform after an $a$-shuffle for a 
	deck of distinct cards and 3 methods of dealing bridge; $a=16$ 
	means 4 riffle suffles, $a=32$ is 5 riffle shuffles, etc.}
\end{table}

\item {\bf How good is the GSR model?}

The GSR model represents idealized riffle shuffling, in the sense that 
every possible shuffle is equally likely.  Human shufflers vary in their 
skill and ``neatness,'' sometimes clumping cards together too much, 
sometimes not enough.  The question is whether the conclusions drawn 
here from the GSR model will still hold when the model is replaced with 
the way real people shuffle cards.  This is a topic for future work.

\end{enumerate}

\paragraph{Acknowledgements.} The authors would like to thank the 
referees, and also Divakar Viswanath and Jeffrey Lagarias, for many 
helpful conversations.

\bibliography{../ref}

\def\lfhook#1{\setbox0=\hbox{#1}{\ooalign{\hidewidth
  \lower1.5ex\hbox{'}\hidewidth\crcr\unhbox0}}}
\providecommand{\bysame}{\leavevmode\hbox to3em{\hrulefill}\thinspace}
\providecommand{\MR}{\relax\ifhmode\unskip\space\fi MR }
\providecommand{\MRhref}[2]{%
  \href{http://www.ams.org/mathscinet-getitem?mr=#1}{#2}
}
\providecommand{\href}[2]{#2}
\begin{thebibliography}{10}

\bibitem{AldousWalks}
D.~Aldous, Random walks on finite groups and rapidly mixing {M}arkov chains, in
  \emph{Seminar on Probability, XVII}, Lecture Notes in Mathematics, vol. 986,
  Springer, Berlin, 1983, 243--297.

\bibitem{ADSLong}
S.~Assaf, P.~Diaconis, and K.~Soundararajan, A rule of thumb for riffle
  shuffling, \emph{Advances in Applied Probability} (to appear).

\bibitem{BayerDiaconis}
D.~Bayer and P.~Diaconis, Trailing the dovetail shuffle to its lair, \emph{Ann.
  Appl. Probab.} \textbf{2} (1992) 294--313.

\bibitem{Carlitz}
L.~Carlitz, Eulerian numbers and polynomials, \emph{Math. Mag.} \textbf{32}
  (1958/1959) 247--260.

\bibitem{Ciucu}
M.~Ciucu, No-feedback card guessing for dovetail shuffles, \emph{Ann. Appl.
  Probab.} \textbf{8} (1998) 1251--1269.

\bibitem{Thesis}
M.~Conger, Shuffling decks with repeated card values, Ph.D. dissertation,
  University of Michigan, Ann Arbor, MI, 2007.

\bibitem{Shuffle1}
M.~Conger and D.~Viswanath, Riffle shuffles of decks with repeated cards,
  \emph{Ann. Probab.} \textbf{34} (2006) 804--819.

\bibitem{Shuffle3}
\bysame, Shuffling cards for blackjack, bridge, and other card games (2006),
  available at http://arxiv.org/abs/math/0606031v1.

\bibitem{Diaconis88}
P.~Diaconis, \emph{Group Representations in Probability and Statistics},
  Institute of Mathematical Statistics Lecture Notes---Monograph Series, 11,
  Institute of Mathematical Statistics, Hayward, CA, 1988.

\bibitem{Diaconis01}
\bysame, Mathematical developments from the analysis of riffle shuffling, in
  \emph{Groups, Combinatorics and Geometry ({D}urham, 2001)}, A.~A. Ivanov,
  M.~W.~Liebeck, and J.~Saxl, eds., World Scientific, River Edge, NJ, 2003,
  73--97.

\bibitem{FulmanCuts}
J.~Fulman, Affine shuffles, shuffles with cuts, the {W}hitehouse module, and
  patience sorting, \emph{J. Algebra} \textbf{231} (2000) 614--639.

\bibitem{Gilbert}
E.~N. Gilbert, Theory of shuffling, Tech. Report MM-55-114-44, Bell Telephone
  Laboratories, New York, October 21, 1955.

\bibitem{ConcreteMath}
R.~L. Graham, D.~E. Knuth, and O.~Patashnik, \emph{Concrete Mathematics}, 2nd
  ed., Addison-Wesley, Reading, MA, 1994.

\bibitem{Mann1}
B.~Mann, How many times should you shuffle a deck of cards? in \emph{Topics in
  Contemporary Probability and Its Applications}, Probability and Stochastics
  Series, CRC, Boca Raton, FL, 1995, 261--289.

\bibitem{Markov}
A.~A. Markov, Extension of the law of large numbers to dependent events,
  \emph{Bull. Soc. Phys. Math.} \textbf{15} (1906) 135--156.

\bibitem{Poincare}
H.~Poincar{\'{e}}, \emph{Calcul des Probabilit{\'{e}}s}, deuxi{\`{e}}me ed.,
  Gauthier-Villars, Paris, 1912.

\bibitem{Reeds}
J.~Reeds, Unpublished manuscript, 1981.

\bibitem{Tanny}
S.~Tanny, A probabilistic interpretation of {E}ulerian numbers, \emph{Duke
  Math. J.} \textbf{40} (1973) 717--722.

\bibitem{Trefethen}
L.~N. Trefethen and L.~M. Trefethen, How many shuffles to randomize a deck of
  cards? \emph{Proc. Roy. Soc. London Ser. A} \textbf{456} (2000) 2561--2568.

\end{thebibliography}

\bigskip

\noindent{\bf Mark Conger} received his B.A. from Williams College in 
1989 and his Ph.D. from the University of Michigan in 2007.  In between 
he worked as a professional programmer for many years.  He enjoys 
woodworking and taking things apart. He currently teaches at the 
University of Michigan.

\noindent\textit{Department of Mathematics, University of Michigan, 2074 
East Hall, 530 Church Street, Ann Arbor, MI 48109}\\
{\tt mconger@umich.edu}

\bigskip

\noindent{\bf Jason Howald} received his B.A. from Miami University in 
1995 and his Ph.D. from the University of Michigan in 2001.  He enjoys 
juggling, baking, and computer programming. He currently teaches at the 
State University of New York at Potsdam.

\noindent\textit{Department of Mathematics, SUNY Potsdam, 44 Pierpont 
Avenue, Potsdam, NY 13676}\\
{\tt howaldja@potsdam.edu}

\end{document}